\documentclass[12pt]{amsart}
\usepackage{amssymb}
\newtheorem{theorem}{Theorem}[section]
\newtheorem{corollary}[theorem]{Corollary}
\newtheorem{proposition}[theorem]{Proposition}

\newtheorem{lemma}[theorem]{Lemma}

\theoremstyle{remark}
\newtheorem{remark}[theorem]{Remark}

\newcommand\be{\begin{equation}}
\newcommand\ee{\end{equation}}
\newcommand\M{\mathcal{M}}
\renewcommand\L{\mathcal{L}}

\renewcommand{\O}{\mathcal{O}}
\newcommand{\Co}{\mathcal{C}}

\newcommand{\U}{\on{U}}

\newcommand{\R}{\mathbb{R}}
\newcommand{\C}{\mathbb{C}}
\newcommand{\cC}{\mathcal{C}}

\newcommand\lie[1]{\mathfrak{#1}}
\renewcommand{\k}{\lie{k}}

\newcommand{\g}{\lie{g}}

\renewcommand{\t}{\lie{t}}

\newcommand{\on}{\operatorname}
\newcommand{\ad}{\on{ad}}

\newcommand{\Ad}{ \on{Ad} }

 \newcommand{\Gl}{ \on{Gl}}

\newcommand{\D}{ \mathcal{D} }

\newcommand\dirac{/\kern-1.2ex\partial} 
\newcommand\qu{/\kern-.7ex/} 

\newcommand{\labell}\label

\renewcommand\a{\mathfrak{a}}

\newcommand{\hra}{\hookrightarrow}

\renewcommand{\d}{{\mbox{d}}}
\newcommand{\ol}{\overline}

\newcommand\sig{\sigma}
\newcommand\eps{\epsilon}
\newcommand\Om{\Omega}
\newcommand\om{\omega}

\newcommand{\f}{\frac}

\newcommand{\p}{\partial}
\renewcommand{\l}{\langle}
\renewcommand{\r}{\rangle}
\newcommand{\hh}{{\textstyle \f{1}{2}}}
\newcommand{\ts}{\textstyle}
\newcommand{\ti}{\tilde}
\newcommand{\olt}{{{\theta}^\dagger}} 

\newcommand\beqn{\begin{equation}}      
\newcommand\eeqn{\end{equation}}      
\newcommand{\ca}{\mathcal}

\renewcommand{\cC}{\mathcal{C}}

\newcommand{\mf}{\mathfrak}
\newcommand{\n}{\mf{n}}
\newcommand{\beq}{\begin{eqnarray*}}
\newcommand{\eeq}{\end{eqnarray*}}

\begin{document}

\title[Linearization of Poisson moment maps]{Linearization of Poisson
actions\\ and  singular values of matrix products}

\date{December 2000}

\author{A. Alekseev}
\address{Institute for Theoretical Physics \\ Uppsala University \\
Box 803 \\ \mbox{S-75108} Uppsala \\ Sweden}
\email{alekseev@teorfys.uu.se}


\author{E. Meinrenken}
\address{University of Toronto, Department of Mathematics,
100 St George Street, Toronto, Ontario M5S3G3, Canada }
\email{mein@math.toronto.edu}

\author{C. Woodward}
\address{Mathematics-Hill Center, Rutgers University,
110 Frelinghuysen Road, Piscataway NJ 08854-8019, USA}
\email{ctw@math.rutgers.edu}

\maketitle

\begin{abstract}
We prove that the linearization functor from the category of
Hamiltonian $K$-actions with group-valued moment maps in the sense of
Lu, to the category of ordinary Hamiltonian $K$-actions, preserves
products up to symplectic isomorphism.  As an application, we give a
new proof of the Thompson conjecture on singular values of matrix
products and extend this result to the case of real matrices.  We give
a formula for the Liouville volume of these spaces and obtain from it
a hyperbolic version of the Duflo isomorphism.
\end{abstract}

\section{Introduction}
\label{sec:intro}
Poisson-Lie groups were introduced by Drinfeld \cite{dr:qu} as
semiclassical analogs of quantum groups. By definition, a Poisson-Lie
group is a Lie group endowed with a Poisson structure such that group
multiplication is a Poisson map.  Poisson-Lie
groups have been used to generalize the Kostant nonlinear convexity
theorem \cite{ko:on,lu:no,fl:pc}, explain the properties of Kostant
harmonic forms on flag manifolds \cite{ko:li,ev:po}, and understand
the symmetries of certain integrable systems \cite{se:dr}.  An
important role in these applications is played by the notion of a {\em
moment map} for a Poisson action of a Poisson-Lie group, due to
J.-H. Lu \cite{lu:mo}. In contrast to ordinary moment maps taking
values in the dual of the Lie algebra, moment maps in the sense of Lu
take values in the {\em dual} Poisson-Lie group.

Compact Lie groups $K$ carry a distinguished non-trivial Lie-Poisson
structure known as the Lu-Weinstein \cite{lu:po} Poisson structure.
For this case, the first author showed \cite{al:po} 
that the categories of symplectic $K$-manifolds with moment maps in 
the dual group $K^*$,
respectively dual of the Lie algebra $\k^*$ are equivalent.  That is,
for every Poisson $K$-action on a symplectic manifold $(M, \Omega)$
with $K^*$-valued moment map $\Psi$, there is a different symplectic
form $\omega$ for which the action is Hamiltonian in the usual sense,
with a $\k^*$-valued moment map $\Phi$. Poisson reductions of $(M,
\Omega,\Psi)$ are isomorphic to reductions of its {\em linearization}
$(M, \omega,\Phi)$ as (stratified) symplectic spaces.

The categories of symplectic $K$-manifolds with $\k^*$- and
$K^*$-valued moment maps have natural structures of {\em tensor
categories}\,: There are operations of products, sums and conjugation
satisfying the usual axioms.  The first main result of this paper is
that the linearization functor preserves these operations up to
symplectomorphism.  The proof is based on a simple Moser isotopy
argument.  As an application, we prove the Thompson conjecture on
singular values of products of complex matrices, which was first
established in a recent paper by Klyachko \cite{kl:ra}, and also the
corresponding statement for real matrices (Theorem \ref{th:th1}).
Independently, a completely different proof of these
results was obtained by Kapovich-Leeb-Millson \cite{ka:po}.

The second main result is a formula comparing the Liouville volume
forms defined by $\om$ and $\Om$.  This formula involves the modular
function for $K^*$ and a Duflo factor. As a corollary, we obtain
Klyachko's formula \cite{kl:ra} for random walk distributions, which
we interpret as a hyperbolic version of the Duflo theorem. That is, a
certain linear map between spaces of compactly supported distributions
on $\k$ and $K^*$ becomes a ring homomorphism (with respect to
convolution) if restricted to $K$-invariants.

\section{Moment maps for Poisson actions}
\label{sec:moment}

In this Section we recall the theory of moment maps 
for Poisson actions of compact Poisson-Lie groups
on symplectic manifolds developed by Lu \cite{lu:mo}.

\subsection{Poisson-Lie groups}
Recall that a Poisson-Lie group is a Lie group $K$ together with a
Poisson bivector $\pi_K$ such that group multiplication is a Poisson
map. This condition implies that the inversion map $K\to K,\, k\mapsto k^{-1}$
is anti-Poisson. The Poisson bivector $\pi_K$ vanishes at the group
unit of $K$, and its linearization $\delta: \ \k \to \k \otimes \k$ is
a 1-cocycle on $\k$. The dual map $\delta^*$ defines a Lie algebra
structure on $\k^*$. The connected, simply-connected Lie group 
$K^*$ with Lie algebra $\k^*$ is called the Poisson dual of 
$K$. It is a Poisson-Lie group, with Poisson bracket induced by the 
Lie algebra structure on $\k$.  

Let the vector space $\g=\k\oplus\k^*$ be equipped with the symmetric
bilinear form $\langle\cdot,\cdot\rangle$ for which $\k$ and $\k^*$
are isotropic and which extends the natural pairing between elements
in $\k$ and $\k^*$.  According to \cite[Theorem 1.12]{lu:po} there is
a unique Lie algebra structure on $\g = \k \oplus \k^*$ for which $\k$
and $\k^*$ are subalgebras and the pairing $\langle\cdot,\cdot\rangle$
is $\g$-invariant.

A Lie group $G$ with Lie algebra $\g$ is called a {\em double} for the
Poisson-Lie group $K$ if the subalgebras $\k,\k^*\to\g$ exponentiate
to closed subgroups $K,K^*\to G$, and the multiplication map
$K^*\times K\to G,\ (l,k) \mapsto lk$ is a diffeomorphism.  In this
case, the left-action of $G$ on itself induces an action on
$K^*=G/K$. Its restriction $K\times K^*\to K^*,\ (k,{l})\mapsto {l}^k$
is called the dressing action of $K$ on $K^*$.  Similarly, the
right-action of $G$ restricts to the dressing action $K^*\times K\to
K,\ ({l},k)\mapsto k^{l}$ on $K=K^*\backslash G$. The two actions are
related by
\begin{equation} \label{kprime}
k {l} = {l}^k k^{l} .\end{equation}
The classification of Poisson-Lie structures on compact, connected Lie
groups $K$ was carried out by Levendorskii and Soibelman \cite{le:al}.
Besides the trivial structure, there is a distinguished example called
the Lu-Weinstein structure. Let $\g=\k^\C$, viewed as a real Lie
algebra, and $\g=\k\oplus\a\oplus\n$ an Iwasawa decomposition.

For any invariant inner product $B$ on $\k$, with complexification 
$B^\C$, the bilinear form 
$$\l\cdot,\cdot\rangle=2\on{Im}\,B^\C$$ 
defines a non-degenerate pairing between $\k$ and $\a\oplus\n$,
identifying $\k^*\cong\a\oplus\n$.  The induced Lie algebra structure
on $\k^*$ defines the Lu-Weinstein Poisson structure on $K$, with
Poisson dual $K^*=AN$ and double $G=K^\C=KAN$.

\subsection{Poisson actions} 
Let $(K,\pi_K)$ be a connected Poisson-Lie group, with Poisson-dual $K^*$, 
and suppose $K$ admits a double $G=K^*K$.
Denote by $\theta^L,\theta^R\in\Om^1(K^*)\otimes\k^*$ 
the left- and right-invariant Maurer-Cartan forms. 
\cite[Corollary 3.6]{lu:mo} states that for every Poisson map 
$\Psi:\,M\to K^*$ from a Poisson manifold $(M,\pi)$ to $K^*$, 
the formula
\begin{equation} \label{eq:momcon}
\xi_M = \pi^\sharp \Psi^*  \langle \theta^R, \xi \rangle, \ \ \xi\in\k
\end{equation}
defines a Lie algebra action of $\k$ on $M$,
i.e. $[\xi_M,\eta_M]=[\xi,\eta]_M$. If this action integrates to a
$K$-action, with generating vector fields $m\mapsto \f{d}{d
t}|_{t=0}\exp(-t\xi).m$ equal to $\xi_M(m)$, then the the triple
$(M,\pi,\Psi)$ is called a Hamiltonian $K$-space with $K^*$-valued
moment map $\Psi$. It follows from the moment map condition
\eqref{eq:momcon} that the action map $K\times M\to M$ is Poisson
\cite[Corollary 3.6]{lu:mo} and that the moment map is $K$-equivariant
\cite[Theorem 4.8]{lu:mo}.  For $\pi_K=0$ this reduces to the usual
definition of a Hamiltonian $G$-space with $\k^*$-valued moment
map. In the special case where $\pi_K$ is the inverse of a symplectic
structure $\Om\in\Om^2(M)$, the moment map condition is equivalent to
\begin{equation}\label{eq:moment}
 \iota(\xi_M)\Om=\Psi^* \l \theta^R,\xi\r.
\end{equation}

There are sum, product, and conjugation operations for Hamiltonian 
$K$-manifolds with $K^*$-valued moment maps, as follows. 
Sum is given by disjoint union.  The product of two Hamiltonian 
$K$-manifolds with $K^*$-valued moment maps
$(M_1,\pi_1,\Psi_1)$ and $(M_2,\pi_2,\Psi_2)$ is given by 
$$ (M_1 \times M_2,\pi_1 + \pi_2,\Psi_1\Psi_2).$$
Indeed, by Flaschka-Ratiu \cite[Lemma 22.3]{fl:pc} the infinitesimal 
action generated by the Poisson map $\Psi_1\Psi_2$ exponentiates to 
the following $K$-action on $M_1\times M_2$, 
$$ k.(m_1,m_2) = (k.m_1,k^{\Psi_1(m_1)}. m_2).$$
The twist product is associative.  It defines a 
tensor category structure on 
Hamiltonian $K$-manifolds, with morphisms given
by equivariant Poisson isomorphisms preserving the moment map.

\begin{lemma} For any Hamiltonian $K$-manifolds with 
$K^*$-valued moment map $(M, \pi, \Psi)$ the formula
\begin{equation} \label{twist}
 (k,m) \mapsto k^{\Psi(m)^{-1}}.m \end{equation}
defines a Poisson action on $(M, -\pi)$ with moment map $\Psi^{-1}$.
We call $(M, -\pi, \Psi^{-1})$ the conjugate of $(M, \pi, \Psi)$.
\end{lemma}
\begin{proof} 
First, we check that \eqref{twist} defines an action. Let $K_L,K_R$ be
two copies of $K$ acting on $G$ by $(k,g)\mapsto kg$ and $(k,g)\mapsto g
k^{-1}$, respectively. Consider $G$ as a $K_L$-equivariant principal
$K_R$-bundle over $K^*=G/K_R$, and let $\Psi^*G$ denote the pull-back
to $M$. The action of $K_L$ on $\Psi^*G$ is free, and has 
$$ \iota:\ M\to \Psi^*G,\ \ m\mapsto (m,\Psi(m))$$
as a cross-section. Using $\iota$ identify $\Psi^*G/K_L=M$. 
We claim that the induced action of $K_R$ on $M$ is the twisted 
$K$-action.  Given $m\in M$ we compute 
$$ (m,\Psi(m) k^{-1}) = (m, ( k^{\Psi^{-1}(m)} )^{-1} ((\Psi^{-1}(m))^k)^{-1})
$$
The action of $k^{\Psi^{-1}(m)}$ takes this back to
$\iota(M)$, which proves the claim.

Since the inversion
map on $K^*$ is anti-Poisson, $\Psi^{-1}$ is a Poisson map for 
the reversed Poisson structure $-\pi$ on $M$.  
We check it is a moment map for the twisted action. 
Let $\on{pr}_\k:\,\g\to\k$ denote projection 
along $\k^*$. Using the moment map condition for $\Psi$, 
\beq
 -\pi^\sharp (\Psi^{-1})^* \langle \theta^R, \xi \rangle (m) &=&
\pi^\sharp \Psi^* \langle \theta^L, \xi \rangle(m) \\
&=& \pi^\sharp \Psi^*
\langle \theta^R,\on{pr}_\k (\Ad_{\Psi(m)}\xi) \rangle (m) \\
&=&(\on{pr}_\k(\Ad_{\Psi(m)}\xi))_M(m)\\
&=&(\xi^{\Psi(m)^{-1}})_M(m),
\eeq
which are the generating vector fields for the twisted action.
\end{proof}

Symplectic reduction extends to the setting of Hamiltonian Poisson
actions with $K^*$-valued moment maps. Suppose $M$ is symplectic
structure and that the action is proper.  For any ${l} \in K^*$, define
$$ M_{l} = \Psi^{-1}({l})/K_{l} \cong \Psi^{-1}(K{l})/K $$
where $K{l}$ is the orbit of ${l}$ under the dressing action of $K$ on
$K^*$. Then $M_{l}$ is a symplectic manifold, if the action of $K$
on $\Psi^{-1}(K{l})$ is free \cite[Theorem 4.12]{lu:mo}.

\subsection{Anti-Poisson involutions}

Recall the definition of compatible involutions from O'Shea-Sjamaar
\cite{os:th}.  Let $K$ be a connected Lie group, together with an
involutive automorphism $\sig_K$. Let $\sig_\k$ denote the
corresponding Lie algebra involution, and define an involution on
$\k^*$ by $\sig_{\k^*}=-(\sig_\k)^*$.
 
An involution $\sig_M:\,M\to M$ of a symplectic manifold $(M,\om)$ 
is called anti-symplectic if $\sig_M^*\om=-\om$. If $M$ carries 
a Hamiltonian $K$-action, with moment map $\Phi:\,M\to\k^*$, then 
$\sig_M$ is called compatible with $\sig_K$ if 
$$\Phi\circ \sig_M = \sig_{\k^*}\circ \Phi.$$  
As explained in \cite{os:th}, since $K$ is connected this implies 
$$\sig_M(k.m) =\sigma_K(k).\sig_M(m).$$
If $M\subset\k^*$ is a coadjoint orbit with the
Kirillov-Kostant-Souriau symplectic structure, such that $M$ is
invariant under $\sig_{\k^*}$, the involution $\sig_M=\sig_{\k^*}|_M$
is compatible with $\sigma_K$.

Suppose $K$ is compact.  Choose a Cartan subalgebra $\t$ of $\k$ such
that $\k^\sig \cap \t$ has maximal dimension.  Let $\t_+^*\subset\t^*$
be a positive Weyl chamber.  The following theorem of O'Shea-Sjamaar
describes the image of the fixed point manifold $M^\sig$ under the
moment map.  The special case where $K$ is a torus and $\sigma_K(k) =
k^{-1}$ is due to Duistermaat \cite{du:co}.
\begin{theorem}[O'Shea-Sjamaar] \label{th:oshea}
Let $(M,\omega,\Phi)$ be a symplectic Hamiltonian $K$-manifold with
proper $\k^*$-valued moment map, and let $\sig_M$ be a
$\sig_K$-compatible anti-symplectic involution on $M$. Then
$$ \Phi(M^\sig) \cap \t^*_+ = \Phi(M)^\sig\cap \t^*_+.$$
\end{theorem}
A theorem of Kirwan says that if $M$ is compact and connected,
$\Delta(M) = \Phi(M) \cap \t^*_+$ is a convex polytope. By Theorem
\ref{th:oshea}, $\Phi(M^\sig) \cap \t^*_+$ is also a polytope,
obtained from the Kirwan polytope by intersecting with the subspace
$(\k^*)^\sig$.

We generalize these definitions to Poisson actions and $K^*$-valued
moment maps as follows.  Let $K$ be a connected Poisson-Lie group,
together with an anti-Poisson involutive automorphism $\sig_K$.  Then
$\sig_{\k^*}$ is a Lie algebra automorphism, and therefore
exponentiates to a Lie group automorphism $\sig_{K^*}$ on the Poisson
dual $K^*$.  For any Hamiltonian Poisson $K$-manifold $(M,\pi, \Psi)$
we say that an anti-Poisson involution $\sig_M$ of $M$ is compatible
with $\sigma_K$ if
\begin{equation} \label{eq:compat}
\Psi\circ \sig_M = \sigma_{K^*}\circ \Psi
.\end{equation}
Since $K$ is connected, this implies $\sig_M(k.m) =
\sigma_K(k).\sig_M(m)$.  Indeed, for anti-Poisson involutions $\sig_M$
and $\sig_K$, the composition $\sig_{K^*} \circ\Psi\circ \sig_M$ is
Poisson, and is the moment map for the action, $(k,m)\mapsto
\sig_M(\sig_K(k).\sig_M(m))$. Condition \eqref{eq:compat} implies that
these are the original moment map and action.  Examples of Hamiltonian
$K$-spaces with compatible involution are $\sig_{K^*}$-invariant
dressing orbits $M$ for the action of $K$ on $K^*$, with
$\sig_M=\sig_{K^*}|_M$.  If $\sig_M$ is a compatible involution of
$(M,\pi,\Psi)$ then it is also a compatible involution of the
conjugate $(M,-\pi,\Psi^{-1})$. Similarly, if $\sig_{M_j}$ ($j=1,2$)
are compatible involutions of $(M_j,\pi_j,\Psi_j)$, then
$\sig_{M_1}\times\sig_{M_2}$ is a compatible involution of their
product.

The fixed point set $M^\sig$ carries an action of the group
$K^\sigma$. For ${l} \in (K^*)^\sigma$ we denote by $M^\sig_{l}$ the
quotient
\begin{equation} \label{fixquot}
 M^\sig_{l} = \Psi^{-1}({l})^\sig / K_l^\sigma .\end{equation}

\subsection{Examples of anti-Poisson involutions}
\label{subsec:ex}
Suppose $K$ is a compact Lie group, equipped with 
the Lu-Weinstein Lie-Poisson structure  
corresponding to an invariant inner product 
$B$ on $\k$. Let $G=K^\C$. 
\begin{lemma}
Let $\sig_\g$ be an anti-linear 
involutive automorphism of $\g=\k^\C$ preserving 
$\k,\k^*$. Suppose $\sig_\k=\sig_\g|_\k$ is an isometry. Then  
the exponentiated automorphism $\sig_K$ of $K$ is 
an anti-Poisson involution. 
\end{lemma}
\begin{proof} Since $\sigma_\g$ preserves $B$ and is anti-linear
it takes $B^\C$ to its complex conjugate and so changes the sign of
$\l\ , \ \rangle = 2\on{Im}(B^\C)$. It follows that the  
involutions $\sig_\k=\sig_\g|_\k$ and $\sig_{\k^*}=\sig_\g|_{\k^*}$
are related by $\sig_{\k^*}=-(\sig_\k)^*$. Therefore $\sig_\k$ changes 
the sign of the cocycle $\delta$ dual to the bracket on $\k^*$.
\end{proof}
We remark that if $\k$ is simple then any involution  
preserves the Killing form, hence also $B$. 
For $\k$ semi-simple, anti-holomorphic involutions 
$\sigma_\g$ preserving $\k$ and $\k^*$ arise from 
automorphisms of the Dynkin diagram, as follows. Any 
automorphisms of the Dynkin diagram gives rise to an 
automorphism of the root system. Composing with the map $\alpha
\mapsto - \alpha$, we obtain an automorphism mapping the positive
roots to the negative roots. Let $\zeta$ be the corresponding Lie
algebra automorphism of $\k$, and $\kappa:\,\g\to \g$ the Cartan
involution given by complex conjugation for $\g=\k^\C$.  Then
$\sigma_\g =\kappa \circ \zeta^\C$ is an anti-linear involution 
preserving $\k,\k^*$. 

Consider for example the case $G=\on{Sl}(r, \C)$ with $r\geq 3$. The
trivial automorphism of the Dynkin diagram $A_{r-1}$ induces complex
conjugation on $G$, while the unique non-trivial  automorphism induces 
\begin{equation} \label{exot}
\sigma_G(g) = P (g^\dagger)^{-1} P , \end{equation}
where $P$ is the anti-diagonal matrix $P_{ij}=\delta_{i,n+1-j}$. 

%

\section{Linearization}
\label{sec:linearization}

In this Section we recall the notion of linearization
for Lu-Weinstein moment maps, and then prove that
linearization commutes with product and conjugation up
to symplectomorphism. From now on, 
$K$ denotes a compact, connected, 
Lie group with Lu-Weinstein Poisson-Lie structure, $K^* =
AN$ denotes its Poisson dual, and $G=K^\C=KAN$ the double. 

\subsection{Linearization Theorem}
 In \cite{al:po} the first author
constructed a 1-1 correspondence between Hamiltonian $K$-manifolds
with $\k^*$-valued moment maps and with $K^*$-valued moment maps. To 
set up this correspondence one first needs an equivariant map from 
$\k^*$ to $K^*$. Let
$\kappa:\,\g\to\g$ be the Cartan involution given by complex conjugation 
of $\g=\k^\C$, and let $\dagger:\,\g\to\g$ be the anti-involution
$$ \xi^\dagger = -\kappa(\xi).$$
We also denote by $\dagger$ the induced anti-involution of $G$,
considered as a real group.  For $K=\U(r)$ and $G = \Gl(r,\C)$,
$g^\dagger=\ol{g}^t$.  Let $B^\sharp: \k^* \rightarrow \k$ be the
isomorphism given by $B$. For any $\mu \in \k^*$, the element $g =
\exp(i B^\sharp(\mu)) \in G$ admits a unique 
decomposition $g = {l}{l}^\dagger$, for some ${l}
\in K^*$.  It follows from the Iwasawa decomposition that the map
$$ E: \ \k^* \to K^*, \ \ \mu \to {l}$$
is a diffeomorphism. It is $K$-equivariant with respect to the coadjoint
action on $\k^*$ and the left dressing action on $K^*$.

Next, we define a certain 1-form on $\k^*$. 
Recall that $\theta^L\in \Om^1(K^*)\otimes\k^*$ is the left-invariant 
Maurer-Cartan
form, and let $\olt^L$ be its image under the map 
$\dagger:\,\k^*\subset\g \to\g$. Then $B^\C(\theta^L,\olt^L)\in \Om^2(K^*)$ is 
imaginary-valued, and we can define a real-valued 1-form on $\k^*$ by 
\begin{equation}\label{eq:beta}
\beta=\f{1}{2i} \mathcal{H}\big( E^* B^\C(\theta^L,\olt^L) \big)
\in\Om^1(\k^*)
\end{equation}
where $\mathcal{H}:\,\Om^\star(\k^*)\to \Om^{\star-1}(\k^*)$
is the standard homotopy operator for the de Rham differential. 
%
%

%
\begin{proposition}\label{prop:beta}
The 1-form $\beta$ has the following property: 
\begin{equation}\label{eq:contractions1}
\iota(\xi_{\k^*}) \d\beta = E^* \l\theta^R,\xi\r 
-\d\l\cdot,\xi\r,\ \ \xi\in\k
\end{equation}
\end{proposition}
A proof of this Proposition will be given in the appendix. 
Suppose now that $(M,\Omega,\Psi)$ is a Hamiltonian $K$-space with 
$K^*$-valued moment map. Let 
\begin{equation}\label{eq:linear}
\Phi = E^{-1} \circ \Psi ,\ \ \om=\Om-\d\Phi^*\beta.
\end{equation}
As an immediate consequence of Proposition \ref{prop:beta}, the moment
map condition \eqref{eq:moment} for $\Psi$ is equivalent to the moment
map condition $\d\l\Phi,\xi\r=\iota(\xi_M)\om$ for the closed 2-form
$\om$.
\begin{theorem}[Linearization Theorem \cite{al:po}]
\label{th:linearization}
Suppose $M$ is a $K$-manifold. Let $\Om,\om\in\Om^2(M)$ be 
two-forms  and $\Psi:\,M\to K^*,\ \Phi:M\to \k^*$ maps 
related by \eqref{eq:linear}. 
Then $(M,\Om,\Psi)$ is a Hamiltonian $K$-space with $K^*$-valued
moment map if and only if $(M,\om,\Phi)$ is a Hamiltonian $K$-space 
with $\k^*$-valued moment map. 
\end{theorem}
We call $(M,\om,\Phi)$ the linearization of $(M,\Om,\Psi)$. 
For example, linearization of a dressing orbit $\D\subset K^*$ gives
the corresponding co-adjoint orbit $\O=E^{-1}(\D)\subset \k^*$. 
Note also that since the pull-backs of $\Om$ and $\om$ to any 
level surface $\Phi^{-1}(\mu)=\Psi^{-1}({l})$ agree, for $\mu=E({l})$, 
there is a canonical isomorphism of symplectic quotients 
$$ M_\mu\cong M_{{l}}$$
of $(M,\om,\Phi)$  at $\mu$ and of $(M,\Om,\Psi)$ at ${l}$.

\subsection{Linearization commutes with products and conjugation}

Now we describe the interaction of linearization with the product and
conjugation operations. We will need the following Moser isotopy lemma.
\begin{lemma} \label{lem:moser}
Let $(M,\om^s,\Phi^s)$ be a family of compact Hamiltonian
$K$-manifolds, $s\in[0,1]$. For $\xi\in\k$ let $\xi_M^s$ denote the
Hamiltonian vector field for $(M,\om^s,\Phi^s)$.  Suppose $\om^s$ and
$\Phi^s$ depend smoothly on $s$ and that there exists a smooth family
of 1-forms $\alpha^s$ such that
\begin{equation}\label{eq:11}
\dot{\om}^s=\d \alpha^s,
\end{equation}
where the dot stands for $\f{d}{d s}$. Assume that for all 
elements $\xi\in\k^K$, 
\begin{equation}\label{eq:12}
\l\dot{\Phi}^s,\xi\r+\iota(\xi_M^s)\alpha^s=0.
\end{equation}
Then there is a smooth isotopy $\phi^s:\,M\to M$ which intertwines the 
$K$-actions for the parameters $0,s$ and which satisfies
$$ (\phi^s)^*\om^s=\om^0,\ \ \ (\phi^s)^*\Phi^s=\Phi^0.$$
Given a family of anti-symplectic 
involutions $\sig_M^s$ of $(M,\om^s,\Phi^s)$, such that each   
$\alpha^s$ is $\sig_M^s$-anti-invariant, one can arrange that 
$\phi^s\circ \sig_M^0=\sig_M^s\circ \phi^s$.
\end{lemma}

\begin{proof}
For each $s\in[0,1]$ let $j^s:\,M\to \ti{M}:=[0,1]\times M$ be the 
inclusion $j^s(m)= (s,m)$. Equip $\ti{M}$ with the $K$-action such 
that the maps $j^s$ are equivariant, with respect to the 
$K$-action on $M$ defined by $\om^s,\Phi^s$. 
Define $\Phi\in C^\infty(\ti{M})\otimes \k^*$ 
by $(j^s)^*\Phi=\Phi^s$, and let
$$ \ti{\om}=\om+\d s\wedge\alpha\in\Om^2(\ti{M})$$
where $\om,\alpha$ pull-back to $\om^s,\alpha^s$ under $j^s$ and vanish 
on $\f{\p}{\p s}$.  Then \eqref{eq:11} is equivalent to 
\begin{equation}\label{eq:14}
\d\ti{\om}=0
\end{equation}
and \eqref{eq:12} is equivalent to the moment map condition
\begin{equation}\label{eq:12a}
\d\l\Phi,\xi\r=\iota(\xi_{\ti{M}})\ti{\om},\ \xi\in\k^K.
\end{equation}
These two equations also hold for the average of $\ti{\om}$ 
under the $K$-action. Since
$$ L(\xi_{\ti{M}})\om=\iota(\xi_{\ti{M}})\d\om+\d\iota(\xi_{\ti{M}})\om
=-\d s\wedge \iota(\xi_{\ti{M}})\dot{\om},$$
the averaging process changes only $\alpha$, but not $\om$. We may therefore 
assume that $\ti{\om}$ is $K$-invariant.

Let $\ti{X}$ be the unique vector field on $\ti{M}$ such that 
$\iota(\ti{X})\ti{\om}=0$ and $\iota(\ti{X})\d s=1$. It is $K$-invariant,
preserves $\ti{\om}$, and its flow $\ti{\phi}^s$ 
takes the slice at $0$ to that at $s$. Let $\phi^s$ 
be the isotopy of $M$ defined by $\ti{\phi}^s\circ j^0=j^s\circ \phi^s$. 
Then $(\ti{\phi}^s)^*\ti{\om}=\ti{\om}$ implies  
$ (\phi^s)^*\om^s=\om^0$. 
Similarly, for $\xi\in \k^K$ we have 
$$L(\ti{X})\l\Phi,\xi\r=\iota(\ti{X})\d\l\Phi,\xi\r 
=\iota(\ti{X})\iota(\xi_{\ti{M}})\ti{\om}
=0.$$
This shows $(\ti{\phi}^s)^*\l\Phi,\xi\r=\l\Phi,\xi\r$, or equivalently 
$(\phi^s)^*\l\Phi^s,\xi\r=\l\Phi^0,\xi\r$. 
$K$-equivariance of the flow 
$\ti{\phi}^s$ implies that $\phi^s$ intertwines the $K$-actions on $M$ for 
the parameters $0,s$. Since the moment maps are determined up to a 
constant in $(\k^*)^K$, this proves 
$(\phi^s)^*\Phi^s=\Phi^0$. 

In the presence of a family of anti-symplectic involutions with 
$(\sig_M^s)^*\alpha^s=-\alpha^s$, the 2-form $\ti{\om}$ changes sign under
the corresponding involution $\sig_{\ti{M}}$ 
of $\ti{M}$. The vector field $\ti{X}$, and therefore its flow, are 
$\sig_{\ti{M}}$-invariant. Equivalently, $\phi^s\circ \sig_M^0
=\sig_M^s\circ \phi^s$. 
\end{proof}

\begin{theorem}[Linearization commutes with products]\label{th:products}
Let $(M_j,\Om_j,\Psi_j)$ be two compact 
Hamiltonian $K$-spaces with $K^*$-valued 
moment maps and $(M_j,\om_j,\Phi_j)$ their linearizations. Consider 
the products 

\beq (M,\Om,\Psi)&=&(M_1\times M_2,\Om_1+\Om_2,\Psi_1\Psi_2)\\
     (M,\om,\Phi)&=&(M_1\times M_2,\om_1+\om_2,\Phi_1+\Phi_2).
\eeq

The Hamiltonian $K$-space $(M,\om,\Phi)$ is equivariantly symplectomorphic 
to the linearization of $(M,\Om,\Psi)$. 
That is, there exists a diffeomorphism $\phi$ of $M$ which 
takes the diagonal $K$-action to the twisted diagonal action, and satisfies
$$ \phi^*\Om=\om+\d\Phi^*\beta,\ \ \ \phi^*\Psi=E\circ \Phi.$$
\end{theorem}

In particular, this implies that $M_1 \times M_2$ is isomorphic as a
Hamiltonian Poisson manifold to $M_2 \times M_1$, which is not at all
obvious from the definition.  It would be interesting to know whether
the category of Hamiltonian Poisson manifolds admits the structure of
a {\em braided} tensor category.

\begin{proof}
Recall that the definition of a $K^*$-valued moment map depends 
on the inner product $B$ on $\k$.
For any $s>0$ consider the rescaled inner product $B^s=s^{-1} B$, and
let $\zeta^s:\,\k^*\to \k^*,\ \mu\mapsto s\mu$.  Replacing $B$ with
$B^s$ replaces the map $E$ by $E^s=(\zeta^s)^*E$ and 
the form $\beta$ by $\beta^s=s^{-1} (\zeta^s)^*\beta.$ We obtain a 
family $(M_j,\Om^s_j,\Psi^s_j)$ of 
Hamiltonian $K$-spaces with $K^*$-valued moment map (relative to
$B^s$), with 
$$\Om^s_j=\om_j+\d\Phi_j^*\beta^s,\ \ 
\Psi^s_j=E^s\circ \Phi_j.$$
Taking the linearizations of their products 
$$(M,\Om^s,\Psi^s)=
(M_1\times M_2,\Om_1^s+\Om_2^s,\Psi_1^s\Psi_2^s)
$$
we obtain a family of Hamiltonian $K$-spaces $(M,\om^s,\Phi^s)$ where 
\beq E^s\circ \Phi^s&=&(E^s \circ \Phi_1) (E^s \circ \Phi_2) 
\\ \om^s&=&\om+\d(\Phi-\Phi^s)^*\beta^s
  \eeq
Consider the limit $s\searrow 0$. The family of moment
maps $\Phi^s$ extends smoothly to $s=0$ by $\Phi^0=\Phi$.
Since the family of 1-forms $\beta^s$
extends smoothly to $s=0$ by $\beta^0=0$, $\om^s$ extends
smoothly to $s=0$ by $\om^0=\om$. We
thus have a family of compact Hamiltonian $K$-spaces, 
$ (M,\om^s, \Phi^s),\ \ s\in[0,1] $ 
connecting $(M,\om,\Phi)$ with the linearization of $(M,\Om,\Psi)$.
The proof is completed by an application of Lemma \ref{lem:moser}, 
with 
$$\alpha^s =\f{\d}{\d s}(\Phi-\Phi^s)^*\beta^s.$$
To check the condition \eqref{eq:12} for $\xi\in \k^K$, we 
note that the first term vanishes in our case since 
$\l\Phi^s,\xi\r$ is independent of $s$.
Since $\iota(\xi_M^s)(\Phi^s)^*\d\beta=0$, 
also $\xi_M^s$ is independent of $s$, 
and therefore
$$\iota(\xi_M^s)\alpha^s=
\f{\d}{\d s}\iota(\xi_M^s)(\Phi-\Phi^s)^*\beta^s=0.$$
\end{proof}

The following corollary of Theorem \ref{th:products} 
is important for the
proof of the Thompson conjecture in the next Section.

\begin{corollary} \label{cor:key}
Under conditions of Theorem \ref{th:products} the reduced spaces 
$(M_1\times M_2)_l$ at $l\in K^*$
and of $(M_1\times M_2)_\mu$ at $\mu=E^{-1}(l)$ 
are symplectomorphic. 
\end{corollary}

\begin{theorem}[Linearization commutes with conjugation]  
\label{th:conj}
Let $(M,\Om,\Psi)$ be a compact Hamiltonian $K$-manifold 
with $K^*$-valued
moment map, and $(M,\om,\Phi)$ its linearization. Consider the 
conjugates
\beq (M,\ti{\Om},\ti{\Psi})&=&(M,-\Om,\Psi^{-1}),\\
(M,\ti{\om},\ti{\Phi})&=&(M,-\om,-\Phi).
\eeq
There exists an equivariant symplectomorphism between $(M,\ti{\om},\ti{\Phi})$ 
and the linearization of $(M,\ti{\Om},\ti{\Psi})$. That is, there exists a
diffeomorphism $\phi$ of $M$, which intertwines the twisted action with
the original $K$-action on $M$, and satisfies
\beq
 \phi^*\ti{\Psi}&=&E\circ \ti{\Phi},\\
 \phi^*\ti{\Om}&=&\ti{\om}+ \d \ti{\Phi}^*\beta
\eeq
\end{theorem}

\begin{proof} 
We proceed as in the proof of Theorem \ref{th:linearization}. 
Replacing $B$ with $B^s$ we obtain a family $(M,\Om^s,\Psi^s)$ of Hamiltonian 
$K$-manifolds with $K^*$-valued moment maps (relative to $B^s$) with 
$$ \Psi^s=E^s \circ \Phi,\ \ \Om^s=\om+\d\Phi^* \beta^s.$$
Conjugating and linearizing we obtain a family $(M,\ti{\om}^s,\ti{\Phi}^s)$ 
of Hamiltonian $K$-manifolds with 
$$ E^s\circ \tilde{\Phi}^s=(E^s \circ \Phi)^{-1}$$
and 
$$ \tilde{\omega}^s = -\omega-\d(\Phi+\tilde{\Phi}^s)^* \beta^s 
$$
These families extend smoothly to $s = 0$ by 
$\tilde{\omega}^0= -\om$ and $\tilde{\Phi}^0 = -\Phi$, 
and connect the linearization of $(M,\ti{\Om},\ti{\Psi})$ 
with the space $(M,\ti{\om},\ti{\Phi})$. 
Therefore, the claim again follows from Lemma \ref{lem:moser}.
\end{proof}

\subsection{Linearization and anti-symplectic involutions}
Suppose $\sig_K$ is an involution of $K$ of the type described in 
Section \ref{subsec:ex}. That is, the corresponding Lie algebra 
involution $\sig_\k$ is an isometry with respect to $B$, and 
extends to a $\C$-anti-linear involution $\sig_\g$ preserving $\k^*$. 
Letting $\sig_{K^*}$ be the induced involution of $K^*$, we have 
\begin{equation}\label{eq:sigE}
E\circ \sigma_{\k^*}=\sigma_{K^*}\circ E
\end{equation}
by the calculation, 
$$ \exp(i B^\sharp(\sigma_{\k^*}(\mu)))
= \exp( - i \sigma_\k B^\sharp(\mu)) =
\sigma_G(\exp( i B^\sharp(\mu))  ) = \sigma_{K^*}(l)
\sigma_{K^*}(l)^\dagger.$$
The 1-form $\beta$ defined in \eqref{eq:beta} changes its sign,
\begin{equation}\label{eq:sigbeta}
(\sig_{\k^*})^*\beta=-\beta. \end{equation}
Suppose now that $(M,\Om,\Psi)$ is a Hamiltonian 
$K$-space with $K^*$-valued moment map, and $(M,\om,\Phi)$ its 
linearization. Equations \eqref{eq:sigE} and \eqref{eq:sigbeta} show that 
an involution $\sig_M$ of $M$ is anti-symplectic for $\Om$ 
and $\sig_K$-compatible with $\Psi$ if and only if it 
is anti-symplectic for $\om$ and $\sig_K$-compatible for
$\Phi$. For $\mu\in (\k^*)^\sig,\ l=E^{-1}(\mu)\in (K^*)^\sig$
one has a homeomorphism of quotients, 
$$  M^\sig_\mu\cong M^\sig_l.$$ 
Using the last part of Lemma \ref{lem:moser} one obtains the following extensions of Theorem \ref{th:products}, \ref{th:linearization}.
In Theorem \ref{th:products}, given 
$\sigma_K$-compatible anti-symplectic involutions
$\sig_{M_1},\sig_{M_2}$, the diffeomorphism $\phi$ can be chosen to be
$\sig_M=(\sig_{M_1},\sig_{M_2})$-equivariant, and  
assuming $\mu\in(\k^*)^\sig$, one has a homeomorphism 
\begin{equation}\label{eq:key}
(M_1\times M_2)^\sig_{l}\cong (M_1\times M_2)^\sig_{\mu}.
\end{equation}
Similarly, in Theorem \ref{th:conj} the diffeomorphism 
$\phi$ can be chosen to be equivariant with respect to a given 
anti-symplectic involution $\sig_M$.

\section{The Thompson conjecture for complex and real matrices}
\label{sec:singular}

In this Section we apply our results to give a new proof of the
Thompson conjecture on singular values of complex matrices and to
extend this result to real matrices. 

\subsection{Moduli spaces for additive and multiplicative problems}
Let $\O_1,\ldots,\O_n\subset \k^*$
be given coadjoint orbits, and $\D_j=E(\O_j)\subset K^*$ 
the corresponding dressing orbits. Also let 
$\cC_i = \D_i K=Kg_iK \subset G$ denote the double coset 
containing $\D_i$. Consider the following three moduli spaces, 
\beq 
\M_\O&=&\{(\xi_1,\ldots,\xi_n)\in \O_1\times\ldots\times \O_n|\ 
\xi_1+\cdots+\xi_n=0\}/K,\\
\M_\D&=&\{(g_1,\ldots ,g_n)\in \D_1\times\ldots\times \D_n|\ 
g_1 \cdots g_n =e\}/K,\\
\M_\Co&=&\{(g_1,\ldots ,g_n)\in \Co_1\times\ldots\times \Co_n|
g_1 \cdots g_n =e\}/K^n,
\eeq
where in the last line $K^n$ acts as follows, 
$$ (k_1,\ldots,k_n).(g_1,\ldots,g_n)= 
(k_1 g_1 k_2^{-1},k_2 g_2 k_3^{-1},\ldots,k_ng_nk_1^{-1}).$$
\begin{lemma}
The natural map $\M_\D\to \M_\Co$ is a homeomorphism. 
\end{lemma}

\begin{proof}
Given a solution $(g_1,\ldots,g_n)\in \Co_1\times\ldots\times \Co_n$ 
with product $\prod g_j=e$, define $k_j\in K$ recursively as follows:
put $k_1=e$, let $k_2\in K$ be the
unique element with $g_1 k_2^{-1}\in K^*$, then let $k_3\in K$ 
the unique element with $k_2 g_2 k_3^{-1}\in K^*$, and so on.
Let $(l_1,\ldots,l_n)\in G^n$ be the image of $(g_1,\ldots,g_n)$ by 
the action of $(k_1,\ldots,k_n)$. By construction $l_j\in K^*$ 
for $j<n$, and since the product is $e$ we must have $l_n\in K^*$. 
This shows that the map $\M_\D\to \M_\Co$ is surjective. 
Starting the recursion with $k_1=k$ rather than $k_1=e$ replaces 
$(l_1,\ldots,l_n)$ by its image under the diagonal dressing action of 
$k$. This shows that the map is a bijection. 
\end{proof}

Corollary \ref{cor:key} shows that there exists a symplectomorphism 
between $\M_\O$ and $\M_\D$. It follows that the three moduli spaces
are all homeomorphic: 
\begin{equation}\label{iso1}
 \M_\O\cong \M_\D\cong \M_\Co.
\end{equation}
Given a $\C$-antilinear involution $\sig_\g$ of $\g=\k^\C$ preserving
$\k,\k^*$ and the inner product $B$ on $\k$, we can similarly 
consider moduli spaces 
\beq 
\M_\O^\sig&=&\{(\xi_1,\ldots,\xi_n)\in \O_1^\sig\times\ldots\times \O_n^\sig|\ 
\xi_1+\cdots+\xi_n=0\}/K^\sig,\\
\M_\D^\sig&=&\{(g_1,\ldots ,g_n)\in \D_1^\sig\times\ldots\times \D_n^\sig|\ 
g_1 \cdots g_n =e\}/K^\sig,\\
\M_\Co^\sig&=&\{(g_1,\ldots ,g_n)\in \Co_1^\sig\times\ldots\times \Co_n^\sig|
g_1 \cdots g_n =e\}/(K^\sig)^n.
\eeq
Again we find $\M_\D^\sig\cong \M_\Co^\sig$, and together with 
\eqref{eq:key}  we obtain homeomorphisms 
\begin{equation}\label{iso2}
\M_\O^\sig\cong \M_\D^\sig\cong \M_\Co^\sig.
\end{equation}

\subsection{Thompson conjecture}
We now specialize to the case of $K=\U(r)$, $G=K^\C=\on{Gl}(r,\C)$.
The Lie algebra $\k$ consists of anti-Hermitian matrices.  Identify
$\k^*$ with Hermitian matrices by the pairing,
$$ \l\mu,\xi\r=\f{1}{i}\on{tr}(\mu\xi).$$
The orbits $\O_j\subset\k^*$ consist of Hermitian matrices with
prescribed eigenvalues $\lambda_j^1,\ldots, \lambda_j^r$.  On the
other hand, the double coset spaces $\Co_j\subset G$ consist of
matrices with positive determinant and prescribed {\em singular}
values $\Lambda_j^1,\cdots, \Lambda_j^r$. (Recall that the singular
values of a matrix $A$ are the eigenvalues of $A A^\dagger$.)
Therefore, the equality of moduli spaces $\M_\O\cong \M_\Co$ has the
following consequence.
\begin{theorem}\label{th:th1}
Let $\lambda_j^k\in \R$, be given real
numbers, $1\le j\le n,\ 1\le k\le r$.  The following four 
conditions are equivalent:
\begin{enumerate}
\item[(a)] there exist complex matrices $A_j$ with singular values
$\exp(\lambda_j^k)$ and product $A_1\cdots A_n=I$;
\item[(b)] there exist self-adjoint matrices $B_j$ with eigenvalues 
$\lambda_j^k$ and sum $ B_1+\cdots + B_n=0;$
\item[(c)] there exist real matrices $A_j$ with singular values
$\exp(\lambda_j^k)$ and product $ A_1\cdots A_n=I$;
\item[(d)] there exist real symmetric matrices $B_j$ with eigenvalues
$\lambda_j^k$ and sum $ B_1+\cdots + B_n=0.$
\end{enumerate}
\end{theorem}

\begin{proof}  The equivalence of (a) and (b), first proved by Klyachko in
\cite{kl:ra}, follows from \eqref{iso1}.  The equivalence of (c) and
(d) follows from \eqref{iso2}.  The equivalence of (b) and (d) follows
from Theorem \ref{th:oshea}, since $\sigma$ acts trivially on the
Cartan in this case.  It was proved independently by Fulton
\cite{fu:ei}.
\end{proof}

We note that in a different work Klyachko \cite{kl:st} gave an
inequality description of the set of coadjoint orbits for which the
additive problem admits a solution. This result was generalized to
arbitrary compact Lie groups by Berenstein-Sjamaar \cite{be:coa}.
Theorem \ref{th:th1} implies the same inequality description for the
multiplicative problem for real matrices.

The more general involutions $\sig_K$ discussed in Section \ref{subsec:ex}
yield ``twisted'' versions of the
Thompson conjecture. For example, from the involution \eqref{exot}
we obtain
\begin{theorem}\label{th:th3}
Let $P$ be the anti-diagonal $n\times n$-matrix $P_{ij}=\delta_{i,n+1-j}$. 
Let $\lambda_j^k\in \R$, $j=1,\ldots\!,n, \ k=1,\ldots,r$ be given real
numbers. Then the following two conditions are equivalent:
\begin{enumerate}
\item[(a)] there exist complex matrices $A_j$ satisfying 
$PA_j^\dagger P=A_j^{-1}$, with 
singular values $\exp(\lambda_j^k)$ and product  
$A_1\cdots A_n=I$;
\item[(b)] there exist self-adjoint matrices $B_j$ anti-commuting with 
$P$, with eigenvalues 
$\lambda_j^k$ and sum $B_1+\cdots + B_n=0$. 
\end{enumerate}
 \end{theorem}
%

Inequality descriptions for additive problems involving involutions 
$\sig_K$ are provided by O'Shea-Sjamaar \cite{os:th}.

\section{Volume forms}
\label{sec:volume}
Let $(M,\Om,\Psi)$ be a Hamiltonian $K$-space with $K^*$-valued moment
map. Since the symplectic form is not preserved by the $K$-action, the
symplectic volume form $(\exp \Om)_{[top]}$ is not $K$-invariant in
general. We will show in this Section that one obtains a $K$-invariant
volume form if one multiplies by the pull-back of a certain
multiplicative character of $K^*$.  Similar volume forms were studied
by Lu in the context of Bruhat-Poisson structures on flag manifolds
\cite{lu:co}.  In the case of dressing orbits, the volumes agree with
the ones considered by Klyachko.

Let $\delta:\ \k\to\wedge^2\k$ be the co-bracket defining the Lu-Weinstein 
structure on $K^*$. It is a 1-cocycle for the adjoint representation 
of $\k$:
\begin{equation}
\label{eq:cocycle}
[\xi, \delta(\eta)]-[\eta, \delta(\xi)] 
-\delta([\xi, \eta])=0,
\end{equation}
using the Schouten bracket on $\wedge\k$. 
%
%
%
The cocycle property \eqref{eq:cocycle} of $\delta$ implies that the operators
$$
\mathbf{L}_\xi = L(\xi_M) - {\ts \f{1}{2}}\iota(\delta(\xi)_M),\ \ 
\xi\in\k
$$
define a representation of $\k$ on the space $\Omega(M)$ of
differential forms. We will construct a differential form $\Gamma$ on
$M$ which is invariant under this $\k$-representation and such that
the top degree part $\Gamma_{[top]}$ is a volume form. Since the
operators $\iota(\delta(\xi)_M)$ lower the degree, $\Gamma_{[top]}$ is
then invariant under the usual $K$-action.

The definition involves the modular function $\tau:\ K^*\to \R_{>0}$
for the group $K^*=AN$, i.e. $\tau(g)$ is the determinant of the 
adjoint representation of $K^*$ on $\k^*$. One finds
$$ \tau(\exp \mu)=e^{-4\pi\l\mu,\rho^\sharp\r},\ \ 
\mu\in {\k}^*.$$
Here $\rho\in\t^*$ is the half-sum of positive roots, and 
$\rho^\sharp=B^\sharp(\rho)\in\t$.

\begin{theorem}
\label{th:invariance}
Let $(M,\Om,\Psi)$ be a $K^*$-valued Hamiltonian $K$-space. 
The differential form $\Gamma=\f{\exp(\Om)}{ \Psi^*\tau^{1/2}  }$ is 
invariant under the action of the operators ${\bf{L}}_\xi$. 
Hence its top form degree part 
$$\Gamma_{[top]}=\f{\exp(\Om)_{[top]}}{\Psi^*\tau^{1/2}}$$
is a $K$-invariant volume form on $M$. 
\end{theorem}
\begin{proof}
The exterior differential of $\tau$ is given by 
$$ \d\tau=-\tau\, 4\pi \,\l\theta^R,\rho^\sharp\r.$$
This shows $L(\xi_{K^*})\tau=-\tau\  
4\pi\,\iota(\xi_M) \l\theta^R, \rho^\sharp\r,$
and together with $\L(\xi_M)\Om=\d\l\Psi^*\theta^R,\xi\r$ 
yields, 
\begin{equation}\label{eq:eqA}
\L(\xi_M)\Gamma=\Gamma\ \Psi^*\big( \d\l\theta^R,\xi\r
+2\pi\, \iota(\xi_{K^*})\l\theta^R,\rho^\sharp\r\big).
\end{equation}
To compute $\iota(\delta(\xi)_M)\Gamma$, observe first that 
by the moment map condition, the contraction of 
$\exp(\Om)$ with any bivector field of the 
form $(\xi_1\wedge\xi_2)_M$ for $\xi_j\in\k$ is given by 
$$ 
\iota((\xi_1\wedge\xi_2)_M)\exp(\Om)
=\Big(\l\Psi^*\theta^R,\xi_1\r \l\Psi^*\theta^R,\xi_2\r
+\iota((\xi_1\wedge\xi_2)_M)\Om\Big)\exp(\Om).
$$
The bivector field $\delta(\xi)_M$ is a linear combination 
of such terms. Using the defining property $\l\mu_1\wedge\mu_2,\delta(\xi)\r
=\l[\mu_1,\mu_2],\xi\r$ of the cocycle, 
the first summand simplifies and we obtain
\begin{equation}\label{eq:eqB}
{\ts \f{1}{2}} \iota(\delta(\xi)_M)\Gamma
=
\Big(\hh \Psi^*\l[\theta^R,\theta^R],\xi\r+
{\ts \f{1}{2}} \iota(\delta(\xi)_M)\Om
\Big)\Gamma
\end{equation}
By the structure equation $\d\theta^R=\hh [\theta^R,\theta^R]$, 
the first terms in \eqref{eq:eqA} and \eqref{eq:eqB} agree. By 
the following Proposition \ref{prop:hard} the second terms agree as well. 
\end{proof}
\begin{proposition}\label{prop:hard}
Let $(M,\Om,\Psi)$ be a $K^*$-valued Hamiltonian $K$-space.  For all
$\xi\in\k$, the contractions of $\Om$ with the bivector field
$\delta(\xi)_M$ are given by the formula, 
$$ \iota(\delta(\xi)_M)\Om=
4\pi\, \iota(\xi_M)\Psi^*\l\theta^R,\rho^\sharp\r.
$$ 
\end{proposition}
The proof of this Proposition is deferred to Appendix \ref{app:prop}.
Now let $J_h^{1/2}:\,\k\to\R_{>0}$ be the unique 
$K$-invariant function 
\begin{equation}\label{eq:hduf} 
J_h^{1/2}(\xi)=\prod_{\alpha\succ 0} 
\f{\on{sinh}\pi \l\alpha,\xi\r}{\pi \l\alpha,\xi\r}\ \ 
\end{equation}
for $\xi\in\t$, where the product is over positive (real) roots of $K$. 
Recall that the Duflo factor $J^{1/2}:\,\k\to\R$
(square root of the Jacobian of the exponential map) is given 
by a similar equation but with $\sin$ rather than $\on{sinh}$. 
We therefore call $J_h^{1/2}$ the {\em hyperbolic Duflo factor}. 
Using the isomorphism $B^\sharp:\,\k^*\to \k$ we will view $J_h^{1/2}$ as 
a function on $\k^*$. 
\begin{theorem}
\label{th:volumelinearization}
Let $(M,\Om,\Psi)$ be a $K^*$-valued Hamiltonian $K$-space, 
and $(M,\om,\Phi)$ its linearization. The top form degree parts of
$\Gamma=\Psi^*\tau^{-1} \exp(\Om)$ and $\exp(\om)$ are related by the 
hyperbolic Duflo factor: 
$$\Gamma_{[top]}=\Phi^*J_h^{1/2} \exp(\om)_{[top]}$$
\end{theorem}
\begin{proof}
Since both sides are $K$-invariant, it suffices to verify the 
identity at points of $m\in \Phi^{-1}(\t^*)=\Psi^{-1}(A)$. Let 
$\mu=\Phi(m)\in\t^*$ and $g=\Psi(m)\in A$. Then 
$g=E(\mu)=\exp(i\zeta/2)$ where $\zeta=B^\sharp(\mu)\in\t$, 
and we have
$$ \tau(g)^{1/2}=e^{-2\pi \l\mu,\rho^\sharp\r}
= e^{-2\pi \l\rho,\zeta\r}
$$
Let $U\subset K^*$ be a slice at $\mu$ for the coadjoint-action
on $\k^*$.
There is a splitting, 
\begin{equation}\label{eq:split1}
T_\mu\k^*=
T_\mu U\oplus T_\mu(G\cdot\mu)=
T_\mu U\oplus \k_\mu^\perp
\end{equation}
where $\k_\mu^\perp$ (the orthogonal complement of the isotropy 
algebra) is embedded via the generating vector fields. Let 
$Y=\Phi^{-1}(U)$. By the Guillemin-Sternberg 
symplectic cross-section theorem, 
$Y$ is a symplectic submanifold, and the embedding $\k_\mu^\perp
\to T_mM$ given by the generating vector fields defines an 
$\om$-orthogonal splitting
\begin{equation}\label{eq:split2}
T_m M = T_m Y \oplus \k_\mu^\perp,
\end{equation}
where the 2-form on $\k_\mu^\perp$ is given by the
Kirillov-Kostant-Souriau formula,
$$ \om(\xi_1,\xi_2)=\l\mu,[\xi_1,\xi_2]\r. $$ 
Let $e_\alpha\in\n$ be root vectors for the positive roots $\alpha$, 
normalized by $B(e_\alpha,e_{-\alpha})=1$. 
Then $\on{Re}(e_\alpha),\on{Im}(e_\alpha)$ form a basis of 
$\t^\perp$, and $\k_\mu^\perp$ is the subspace corresponding to 
roots with $\l\alpha,\zeta\r \not=0$. By a short calculation, 
$$ \om_\mu(\on{Re}(e_\alpha),\on{Im}(e_\alpha))
=\pi\l\alpha,\zeta\r.
$$
The splitting \eqref{eq:split2} is also $\Om$-orthogonal. The
pull-backs $\Om_Y$ and $\om_Y$ to $Y$ differ by the pull-back 
by $\Phi|_Y$ of the 2-form $\varpi=\d\beta$. Since 
$\on{ker}(\d_m\Phi)\cap T_mY$ is a co-isotropic subspace 
of $T_mY$ and $\om_Y,\Om_Y$ agree on that subspace, it follows
that the top exterior powers of $\om_Y$ and $\Om_Y$ are equal. 
Therefore, 
$$ 
\exp(\Om_m)_{[top]}=\exp(\om_m)_{[top]}
\prod_{\alpha\succ 0,\ \l\alpha,\zeta\r\not=0}
\f{\Om_g(\on{Re}(e_\alpha),\on{Im}(e_\alpha))}
{\om_\mu(\on{Re}(e_\alpha),\on{Im}(e_\alpha))}.
$$
Since 
$$ \Om_g(e_\alpha,e_{-\alpha})=
e^{-\pi\l\alpha,\zeta\r}
\on{sinh}(\pi \l\alpha,\zeta\r)
$$ 
this gives, 
$$ 
\exp(\Om_m)_{[top]}=\exp(\om_m)_{[top]}\
e^{-2\pi\l\rho,\zeta\r}
\prod_{\alpha\succ 0,\ \l\alpha,\zeta\r\not=0}
\f{\on{sinh}(\pi \l\alpha,\zeta\r)}{\pi \l\alpha,\zeta\r}.
$$
as required.
\end{proof}

\begin{remark}
The proof has not actually used non-degeneracy of the 2-forms 
$\om$ resp. $\Om$. Since 
$\Gamma_{[top]}$ is a volume form if and only if $\Om$ is non-degenerate, 
we have re-proved the second half of Theorem \ref{th:linearization}: 
The 2-form $\om$ of the linearization is non-degenerate if and only 
if the 2-form $\Om$ is non-degenerate.
\end{remark}

\section{DH-measures and the hyperbolic Duflo isomorphism}
\label{sec:dh}
In this Section we identify $\k\cong\k^*$ using $B^\sharp$. In particular
$E:\,\k\to K^*=AN$ is the map such that 
$$\exp(i\xi)=E(\xi)E(\xi)^\dagger$$
for $\xi\in\k$. We will think of $E$ as some kind of exponential map, and 
define a {\em hyperbolic Duflo map} 
$${\bf{D}}_h=E_*\circ J^{1/2}_h:\ \ca{E}'(\k)\to \ca{E}'(K^*)$$
in analogy to the usual Duflo map ${\bf{D}}=
\exp_*\circ J^{1/2}:\ \ca{E}'(\k)\to \ca{E}'(K)$. (Here $\ca{E}'(\cdot)$ 
denotes the space of compactly supported distributions.) Recall that 
${\bf{D}}$ is a ring homomorphism if restricted to invariant distributions. 
Using Theorem \ref{th:products} we will show that the same holds true for the 
hyperbolic Duflo map ${\bf{D}}_h$. 

For any compact $\k$-valued Hamiltonian $K$-space $(M,\om,\Phi)$, 
the Duistermaat-Heckman measure is the compactly supported distribution 
on $\k$ given as a push-forward of the Liouville 
measure, 
$$u:=\Phi_*\ |(e^\om)_{[top]}|\in\ca{E}'(\k)^K$$ 
where $\ca{E}'(\cdot )$ denotes the space of compactly supported
distributions.  Similarly, for the corresponding
$K^*$-valued Hamiltonian $K$-space $(M,\Om,\Psi)$ we define a
DH-measure
$$\mf{m}:=\tau^{-1}\,\Psi_*\ |(e^\Om)_{[top]}|\in \ca{E}'(K^*)^K.$$ 
It is an immediate consequence of \ref{th:volumelinearization} that 
the two measures are related by 
\begin{equation}
\label{eq:dhlinearization} 
\mf{m}={\bf{D}}_h(u). 
\end{equation}
Now suppose $(M_j,\Om_j,\Psi_j)$ are two $K^*$-valued Hamiltonian 
$K$-spaces, and $(M_j,\om_j,\Phi_j)$ their linearizations. Let 
$\mf{m}_j,\,u_j$ denote the respective DH-measures, so that 
$\mf{m}_j={\bf{D}}_h(u_j)$. The DH-measures for the product 
$(M_1\times M_2,\Om_1+\Om_2,\Psi_1\Psi_2)$ is the convolution on 
the group $K^*$, 
$$\mf{m}=\mf{m}_1 *_{K^*} \mf{m}_2,$$
while the DH-measure for 
$(M_1\times M_2,\om_1+\om_2,\Phi_1+\Phi_2)$ is a 
convolution on the vector space $\k$,
$$u=u_1 *_{\k} u_2.$$ 
Since products commute with linearizations up to symplectomorphism 
(Theorem \ref{th:products}), we conclude that $\mf{m}={\bf{D}}_h(u)$. Thus  
\begin{equation}\label{eq:conv}
{\bf{D}}_h(u_1) *_{K^*} {\bf{D}}_h(u_2) = {\bf{D}}_h(u_1 *_{\k} u_2)
\end{equation}
for any distributions $u_1,u_2$ given as DH-measures of Hamiltonian
$K$-spaces. In particular, it holds for DH-measures of coadjoint
orbits; this is one of the results of Klyachko \cite{kl:ra}.  Since
linear combinations of delta distributions are dense in the space
$\ca{E}'(\k)$ of compactly supported distributions, linear
combinations of DH-measures of coadjoint orbits are dense in the space
$\ca{E}'(\k)^K$ of invariant compactly supported distributions, by
averaging.  Therefore, \eqref{eq:conv} holds for {\em arbitrary}
elements $u_j\in\ca{E}'(\k)^K$.  This gives the following
reformulation of Klyachko's result.
\begin{theorem}[Hyperbolic Duflo theorem] \label{th:hduflo}
The map 
$$ {\bf{D}_h}=E_*\circ J_h^{1/2}:\, \ca{E}'(\k)\to \ca{E}'(K^*)$$
is a ring isomorphism if
restricted to $K$-invariant distributions.  That is, \eqref{eq:conv}
holds for all $u_1,u_2\in\ca{E}'(\k)^K$. 
\end{theorem}

\begin{appendix}
\section{Proof of Proposition \ref{prop:beta}}
\label{app:beta}
In this section we prove the property \eqref{eq:contractions1}
of the 1-form $\beta$ used in the linearization construction:
\begin{equation}\label{eq:lim}
\iota(\xi_{\k^*}) \d\beta =2 E^*\,\on{Im} \, B^\C(\theta^R,\xi)  
-\d\l\cdot,\xi\r,\ \ \xi\in\k.
\end{equation}
Let 
$$\Upsilon:\ \k^*\to G,\ \mu\mapsto \exp(i B^\sharp(\mu)) $$
so that $\Upsilon=E E^\dagger$. It is straightforward to check 
$$ 6 \d E^* B^\C(\theta^L,\olt^L)
=\Upsilon^* B^\C(\theta^L,[\theta^L,\theta^L]).$$
(Here and for the rest of this Section $\theta^L,\theta^R$ denote the 
Maurer-Cartan forms for the group $G$. This does not conflict with 
our earlier notation, since the Maurer-Cartan forms for $K^*$  
are given simply by pull-back under the inclusion $K^*\hra G$.)
Hence, 
$$ \d\beta=
-\f{i}{12}\mathcal{H}\big(\Upsilon^* B^\C(\theta^L,[\theta^L,\theta^L])\big)
+\f{i}{2} E^* B^\C(\theta^L,\olt^L).$$
Let us denote the first summand by $\varpi_1$ and the second summand 
by $\varpi_2$. The contractions of $\varpi_2$ with generating vector 
fields $\xi_{\k^*}$ for $\xi\in \k$ are calculated in 
\cite[Lemma 10]{al:mom}
$$ 
\iota(\xi_{\k^*})\varpi_2=\f{i}{2} \Upsilon^* B^\C(\theta^L+\theta^R,\xi)
+2 E^* \,\on{Im} \, B^\C(\theta^R,\xi),\ \ \xi\in\k.
$$
To find the contractions of $\iota(\xi_{\k^*})\varpi_1$, we use the 
identity 
$$ \iota(\xi_G) B^\C(\theta^L,[\theta^L,\theta^L])
=- 6 \d B^\C(\theta^L+\theta^R,\xi).$$ 
Since $\mathcal{H}$ anti-commutes with $\iota(\xi_{\k^*})$, this shows 
\beq  
\iota(\xi_{\k^*})\varpi_1&=&-\f{i}{2}\mathcal{H} \Upsilon^*
\big(\d B^\C(\theta^L+\theta^R,\xi)\big)\\
&=& -\f{i}{2} \Upsilon^* B^\C(\theta^L+\theta^R,\xi)
+ \f{i}{2} \d \mathcal{H}\big( \Upsilon^* B^\C(\theta^L+\theta^R,\xi)\big).
\eeq
From the definition of $\Upsilon$ and of the homotopy operator, 
one finds that 
$$ \mathcal{H}\big( \Upsilon^* B^\C(\theta^L,\xi)\big)
=\mathcal{H}\big( \Upsilon^* B^\C(\theta^R,\xi)\big)
=i \l\cdot,\xi\r.
$$
Hence
$$ \iota(\xi_{\k^*}) \varpi_1=\f{i}{2} \Upsilon^* B^\C(\theta^L+\theta^R,\xi)
    - \d \l\cdot,\xi\r$$ 
Summing with the expression for $ \iota(\xi_{\k^*}) \varpi_2$, we obtain
\eqref{eq:lim}.

\section{Proof of Proposition \ref{prop:hard}}
\label{app:prop}

It is convenient to introduce an orthonormal basis $e_a$
of $\k$. Let $\eps^a\in\k^*\cong\a\oplus\n$ be the dual basis.  We
denote the structure constants of $\k$ by $f_{ab}^c$ and those of
$\k^*$ by $F^{ab}_c$. Thus
\begin{equation}\label{eq:commrel}
[e_a,e_b]=f_{ab}^c e_c,\ \ [\eps^a,\eps^b]=F^{ab}_c \eps^c,\ \ 
[e_a,\eps^b]=-f_{ac}^b \eps^c + F^{bc}_a e_c,
\end{equation}
using summation convention.  Let $v_a=(e_a)_{K^*}$ denote the dressing
vector fields, and $(\eps^a)^R$ the right-invariant vector fields on
$K^*$.  Let $S_{ab}\in C^\infty(K^*)$ be defined by
$$ v_a=S_{ab}(\eps^b)^R.$$ 
In terms of the right-invariant Maurer-Cartan forms, 
$ S_{ab}=\iota(v_a)\theta_b^R$. Note that the restriction 
to any dressing orbit $\D\subset K^*$ is given in terms 
of the symplectic form $\Om$ on $\D$ by $S_{ab}|_\D=-\hh \Om(v_a,v_b)$. 
In particular, $S_{ab}$ is anti-symmetric. 

Recall that $\rho^\sharp=B^\sharp(\rho)\in \t$ where $\rho$ is the
half-sum of positive roots, and write $\rho^\sharp=\rho^b e_b$.
\begin{lemma} \label{lem:1}
$ F^{ab}_a=4\pi \rho^b.$ \end{lemma} 
\begin{proof}
For all $\mu=\mu_a \eps^a\in\k^*$, the number 
$F^{ab}_a\mu_b$ is the trace of the operator $-\ad(\mu)$ on 
$\k^*$. For $\mu\in\mf{n}$, the operator $-\ad(\mu)$ is nilpotent
and therefore has zero trace. Suppose $\mu\in\a$, and 
let $\zeta=B^\sharp(\mu)\in\t$. Since the pairing between $\k^*\supset 
\a=i\t$ and $\k$ is given by $2\on{Im}B^\C$, we have $\mu=\f{i}{2} \zeta$. 

On any root space $\on{span}_\C(e_\alpha)\subset\n\subset\k^*$, 
$-\ad(\mu)$ acts as a scalar $-2\pi i \l\alpha, \f{i}{2} \zeta\r
=\pi \l\alpha,\zeta\r$, hence has trace $2\pi \l\alpha,\zeta\r$.
It follows that the trace of $-\ad(\mu)$ on $\k^*$ is, 
$$ F^{ab}_a\mu_b=2\pi \sum_{\alpha}\l\alpha, \zeta\r=
4\pi \l \rho, \zeta\r 
= 4 \pi \l \mu,\rho^\sharp\r=4\pi \rho^b \mu_b.
$$
\end{proof}
\begin{lemma} \label{lem:2}
$ F_b^{ac}S_{ac}+4\pi \rho^a S_{ab}=0.$
\end{lemma}

\begin{proof}
We claim that the statement of the Lemma is equivalent to the 
equation, 
\begin{equation}\label{eq:equivalent}
(\eps^a)^R S_{ab}=0.
\end{equation}
Indeed, using the definition of $S_{ab}$ we have 
$$ (\eps^a)^R S_{ab}=\iota([(\eps^a)^R,v_a])\theta_b^R
+\iota(v_a)\ L((\eps^a)^R)\theta_b^R.
$$
Since $L((\eps^a)^R)\theta_b^R=F^{ac}_b \theta_c^R$ the second term 
is $F^{ac}_b S_{ac}$. To compute the first term, note that  
the dressing vector fields $v_a$, together with {\em minus} the 
right-invariant vector fields $-(\eps^b)^R$, are the generators for 
the $G$-action on $K^*=G/K$. Therefore, using \eqref{eq:commrel}, 
and Lemma \ref{lem:1},
%
%
%
$[(\eps^a)^R,v_a]=F^{ac}_a v_c=4\pi \rho^a v_a$
which identifies the first term with 
$4\pi \rho^a S_{ab}$.  

It remains to show \eqref{eq:equivalent}. This condition is 
equivalent to the 
vanishing of the second order differential operator
$\Delta_{K^*}=v_b (\eps^b)^R + (\eps^a)^R v_a$ on $K^*$, 
because 
$$ \Delta_{K^*}=(\eps^a)^R S_{ab} (\eps^b)^R-S_{ab} (\eps^a)^R (\eps^b)^R
=((\eps^a)^R S_{ab}) (\eps^b)^R.$$
Let $p:\,G\to K^*=G/K$ be the projection.
Then $ p^* \circ \Delta_{K^*}=\Delta_G \circ p^*$
where 
$$ \Delta_G=(e_a)^R (\eps^a)^R + (\eps^a)^R (e_a)^R 
$$
is the Casimir operator on $G$ corresponding to the 
invariant bilinear form $\l\cdot,\cdot\r$.
Since $\Delta_G$ is $\Ad(G)$-invariant, 
we can replace the superscript ``R'' by a superscript ``L''. 
Hence 
$$ \Delta_G=(e_a)^L (\eps^a)^L + (\eps^a)^L (e_a)^L
=2 (\eps^a)^L (e_a)^L + F^{ar}_a (e_r)^L$$
where we have used $f_{ab}^a=0$. The vector fields $(e_a)^L$
generate the right-$K$ action and therefore  
vanish on right-$K$-invariant functions.  
It follows that $p^* \circ \Delta_{K^*}=\Delta_G\circ p^*=0$, 
so that $\Delta_{K^*}=0$. 
\end{proof}

Now let $S=\hh S_{ab}\eps^a\wedge \eps^b$. The cocycle $\delta(\xi)$ is 
given in terms of the basis by 
$\delta(\xi)=\hh F^{ab}_c \xi^c \,e_a\wedge e_b $. 

\begin{lemma} \label{lem:3}
For all $\xi\in\k$, $\l S,\delta(\xi)\r=
2\pi \l S,\,\xi\wedge \rho^\sharp\r$.  
\end{lemma}

\begin{proof}
Using Lemma \ref{lem:2} we compute, 
$$ 
2\pi \l S,\,\xi\wedge \rho^\sharp\r=2\pi S_{cb}\, \rho^b\xi^c
= \hh S_{ab} F^{ab}_c  \xi^c=\l S,\delta(\xi)\r.
$$ 
\end{proof}

Proposition \ref{prop:hard} is now a direct consequence of Lemma 
\ref{lem:3}, together with the moment map condition.

\end{appendix}

\end{document}